\documentclass[12pt]{article}

\newtheorem{theorem}{Theorem}[section]

\newtheorem{lemma}{Lemma}[section]
\newtheorem{corollary}{Corollary}[section]

\usepackage{amsmath}

\usepackage{natbib}
\usepackage{graphicx}
\usepackage{amssymb,amsmath}
\usepackage{enumerate}
\usepackage{url}

\oddsidemargin
-0.2in \topmargin -.60in
\begin{document}

\title{Small-Deviation Inequalities for Sums of Random Matrices}

\author{
Xianjie~Gao,\quad Chao~Zhang\footnote{Corresponding author},\quad Hongwei~Zhang\\
School of Mathematical Sciences, Dalian University of Technology \\
Dalian, Liaoning, 116024, P.R. China\\
\texttt{xianjiegao@foxmail.com; chao.zhang@dlut.edu.cn; hwzhang@dlut.edu.cn}
}
\date{}



\maketitle

\begin{abstract}%
Random matrices have played an important role in many fields including machine learning, quantum information theory and optimization. One of the main research focuses is on the deviation inequalities for eigenvalues of random matrices. Although there are intensive studies on the large-deviation inequalities for random matrices, only a few of works discuss the small-deviation behavior of random matrices. In this paper, we present the small-deviation inequalities for the largest eigenvalues of sums of random matrices. Since the resulting inequalities are independent of the matrix dimension, they are applicable to the high-dimensional and even the infinite-dimensional cases.

\end{abstract}


\section{Introduction}

Random matrices have been widely used in many problems, for example, the compressed sensing \citep{chandrasekaran2012convex}, the high-dimensional data analysis \citep{buhlmann2011statistics}, the matrix approximation \citep{halko2011finding,gittens2016revisiting} and the dimension reduction \citep{dimensionreduction}. In the literature, one of main research issues is to study the deviation behavior of the eigenvalues (or singular values) of random matrices.

In general, there are two types of deviation results studied in probability theory: one is the large-deviation inequality that describes the behavior of the probability ${\mathbb{P}}(|x|>t)$ for large $t$; and the other is the small-deviation (or small-ball) inequality that controls the probability ${\mathbb{P}}(|x|<\epsilon)$ for small $\epsilon$.

The early large-deviation inequalities for sums of random matrices can be dated back to the work of \citet{ahlswede2002strong}. \citet{tropp2012user} improved their results and developed a user-friendly framework to obtain the large-deviation inequalities for sums of random matrices. To overcome the limitation of the matrix-dimension dependence, \citet{hsu2012tail} and \citet{minsker2017some} introduced the concepts of intrinsic dimension and effective dimension to tighten the large-deviation inequalities, respectively. Moreover, \citet{zhang2017lsv} applied a diagonalization method to obtain the dimension-free large-deviation random for largest singular value of sums of random matrices, while it remains a challenge to select the auxiliary matrices and functions. In the scenario of single random matrix, \citet{ledoux2007deviation} studied the largest eigenvalues of Gaussian unitary ensemble matrices and \citet{vershynin2010introduction} studied the singular values of the sub-Gaussian and sub-exponential matrices.

Small-deviation problems were stemmed from some practical applications, {\it e.g.,}  approximation problem \citep{li1999approximation}, Brownian pursuit problems \citep{li2001capture}, quantization problem \citep{dereich2003link} and convex geometry \citep{klartag2007small}. For more details, we refer to the bibliography maintained by \citet{lifshits2006bibliography}. There have been some works on the small-deviation inequalities for the specific types of random matrices. \citet{aubrun2005sharp} obtained the small-deviation inequalities for the largest eigenvalue of a single Gaussian unitary ensemble matrix. \citet{rudelson2010non} presented the small-deviation inequalities for the smallest singular value of the random matrix with independent entries. \citet{volodko2014small} estimated the small-deviation probability of the determinant of the matrix ${\textbf B}{\textbf B}^T$, where ${\textbf B}$ is a $d \times \infty$ random matrix whose entries obey a centered joint Gaussian distribution. To the best of our knowledge, there are few works on the small-deviation inequalities for sums of random matrices.

\subsection{Related Works}

Let $\{{\textbf X}_1,{\textbf X}_2,\cdots,{\textbf X}_K\}\subset \mathbb{C}^{d \times d}$ be a finite sequence of independent random Hermitian matrices. It follows from Markov's inequality that
\begin{equation*}
   {\mathbb{P}}\left\{\lambda_{\max}\Big(\sum_k{\textbf X}_k\Big)\geq t\right\}\leq \inf_{s>0}\left\{ {\rm e}^{-\theta t}\cdot {\rm tr}\big(\,\mathbb{E}\,{\rm e}^{\theta\sum_k{\textbf X}_k}\big)\right\},
\end{equation*}
where $\lambda_{\max}$ denotes the largest eigenvalue. By using Golden-Thompson inequality, \citet{ahlswede2002strong} bounded the trace of the matrix moment generating function (mgf) in the following way:
\begin{align}\label{eq:aw}
{\rm tr}\big(\mathbb{E}\,{\rm e}^{\theta \sum_k{\textbf X}_k}\big)\leq {\rm tr}({\textbf I})\cdot \Big[\prod_k\lambda_{\max}\big(\mathbb{E}{\rm e}^{\theta {\textbf X}_k}\big)\Big]=d \cdot \exp \left(\sum_k\lambda_{\max}\big(\log \mathbb{E}{\rm e}^{\theta {\textbf X}_k}\big)\right),
\end{align}
where ${\rm tr}({\textbf A})$ stands for the trace of the matrix ${\textbf A}$. By applying Lieb's concavity theorem, \citet{tropp2012user} achieved a tighter matrix mgf bound than the above one:
\begin{equation}\label{eq:tropp}
{\rm tr}\big(\mathbb{E}\,{\rm e}^{\theta \sum_k{\textbf X}_k}\big)\leq d\cdot \exp \left(\lambda_{\max}\Big(\sum_k\log \mathbb{E}{\rm e}^{\theta {\textbf X}_k}\Big)\right),
\end{equation}
where ``the eigenvalue of sum of matrices" is smaller than ``the sum of eigenvalues of matrices" in the right-hand side of \eqref{eq:aw}. However, there still remains a shortcoming that the result \eqref{eq:tropp} is dependent with the matrix dimension $d$, and its right-hand side will become loose for high-dimensional matrices.

To overcome the shortcoming, \citet{hsu2012tail} employed the intrinsic dimension
$\frac{{\rm tr}({\textbf X})}{\lambda_{\max}({\textbf X})}$ to replace the ambient dimension $d$ in the case of real symmetric matrices. \citet{minsker2017some} provided a dimension-free version of Bernstein's inequality for sequences of independent random matrices. \citet{zhang2017lsv} introduced a diagonalization method to obtain the tail bounds for LSV of the sum of random matrices. Although their bounds are independent of the matrix dimension and overcame the aforementioned first shortcoming, there still remains a challenge to select the appropriate parameters to obtain the tighter bounds.

There are also some small-deviation results on one single random matrix. \citet{edelman1988eigenvalues} presented the small-deviation behavior of the smallest singular value of a Gaussian matrix:
\begin{equation*}
  \lim _{d\rightarrow \infty}{\mathbb{P}}\left\{s_{\min}({\textbf A})\leq {\frac{\epsilon}{\sqrt{d}}}\right\}= 1-\exp \Big(-\epsilon-\frac{\epsilon^2}{2}\Big),
\end{equation*}
where ${\textbf A}$ is a $d\times d$ random matrix whose entries are independent standard normal random variables.
\citet{rudelson2008littlewood} studied the the small-deviation bound of the smallest singular value of a sub-gaussian matrix:
\begin{equation*}
  {\mathbb{P}}\left\{s_{\min}({\textbf B})\leq \frac{\epsilon}{ \sqrt{d}}\right\}\leq C\cdot\epsilon+c^d,
\end{equation*}
where $C>0$, $c\in (0,1)$ is only depend on the sub-gaussian moment of its entries and ${\textbf B}$ is a $d\times d$ random matrix whose entries are i.i.d. sub-gaussian random variables with zero mean and unit variance. However, to the best of our knowledge, there is few work on the small-deviation inequalities for sums of random matrices.

\subsection{Overview of Main Results}

In this paper, we present the small-deviation inequalities for the largest eigenvalue for sums of independent random Hermitian matrices, that is, the upper bound of
\begin{equation*}
    {\mathbb{P}}\left\{\lambda_{\max}\Big(\sum_k{{\textbf X}_k}\Big)\leq \epsilon\right\}.
  \end{equation*}
In particular, we first present some basic small-deviation results of random matrices. We then obtain several types of small-deviation inequalities for the largest eigenvalue of sums of independent random positive semi-definite (psd) matrices. Different from the large-deviation inequalities for random matrices, the resulting small-deviation inequalities are independent of the matrix dimension $d$ and thus our finding are applicable to the high-dimensional and even infinite-dimensional cases.

The rest of this paper is organized as follows. In Section \ref{sec:pr}, we introduce some useful notations and then give some basic results on small-deviation inequalities for random matrices. The small-deviation results for sums of random psd matrices are presented in Section \ref{sec:psd}. The last section concludes the paper.

\section{Basic Small-Deviation Inequalities for Random Matrices}\label{sec:pr}

In this section, we first introduce the necessary notations and then present some basic small-deviation results of random matrices.

\subsection{Necessary Notations}

Given a Hermitian matrix ${\textbf A}$, denote ${\lambda_{\max}{(\textbf A})}$ and ${\lambda_{\min}{(\textbf A})}$ as the largest and the smallest eigenvalues of ${\textbf A}$, respectively. Denote ${\rm tr}({\textbf A})$ and $\|{\textbf A}\|$ as the trace and the spectral norm of ${\textbf A}$, respectively. Let ${\textbf I}$ be the identity matrix, ${\textbf U}$ be the unitary matrix and ${\textbf U}^*$ stand for the Hermitian adjoint of ${\textbf U}$.

By the spectral mapping theorem, given a real-value function $f:{\mathbb R}\rightarrow {\mathbb R}$, then
\begin{equation*}
  f({\textbf A})={\textbf U}\cdot f({\Lambda})\cdot {\textbf U}^*,
\end{equation*}
where ${\textbf A}={\textbf U}\Lambda {\textbf U}^*$. If $f(a)\leq g(a)$ for $a\in I$ when the eigenvalues of ${\textbf A}$ lie in $I$, then there holds that $f({\textbf A})\preceq g({\textbf A})$.

\subsection{Basic Small-Deviation Inequalities for Random Matrices}

Subsequently, we come up with the small-deviation inequalities for random matrices. First, we consider a small-deviation bound for one single matrix:
\begin{lemma}\label{lem:laplace}
  Let ${\textbf Y}$ be a random Hermitian matrix. Then for any $\epsilon> 0$,
  \begin{equation*}\label{eq:laplace}
    {\mathbb{P}}\left\{\lambda_{\max}({\textbf Y})\leq \epsilon\right\}\leq \inf_{\theta>0}\left\{\frac{1}{d}\cdot {\rm e}^{\theta \epsilon}\cdot {\mathbb{E}}\,{\rm tr}\big(\,{\rm e}^{-\theta {\textbf Y}}\big)\right\}.
  \end{equation*}
\end{lemma}

{\it Proof:} For any $\theta>0$, we have
    \begin{align*}
{\mathbb{P}}\left\{\lambda_{\max}({\textbf Y})\leq \epsilon\right\}&={\mathbb{P}}\{{\rm e}^{-\lambda_{\max}(\theta{\textbf Y})}\geq {\rm e}^{-\theta\epsilon}\}\\
& \leq {\mathbb{E}}\,{\rm e}^{-\lambda_{\max}(\theta{\textbf Y})}\cdot {\rm e}^{\theta\epsilon} \quad \mbox{[by Markov's inequality]}\\
&={\mathbb{E}}\,{\rm e}^{\lambda_{\min}(-\theta{\textbf Y})}\cdot {\rm e}^{\theta\epsilon}\quad \mbox{[since $-\lambda_{\max}({\textbf A})=\lambda_{\min}({-\textbf A})$]} \\
& ={\mathbb{E}}\,\lambda_{\min}({\rm e}^{-\theta{\textbf Y}})\cdot {\rm e}^{\theta\epsilon}\quad \mbox{[by Spectral mapping theorem]}\\
&\leq \frac{1}{d}\cdot {\rm e}^{\theta \epsilon}\cdot {\mathbb{E}}\,{\rm tr}\big(\,{\rm e}^{-\theta {\textbf Y}}\big).\nonumber\\
\end{align*}
 The last inequality holds because the minimum eigenvalue of a positive definite (pd) matrix is dominated by the ${\rm tr}(\cdot)/d$. Since this inequality holds for any $\theta>0$, taking an infimum over $\theta>0$ completes the proof.\hfill$\blacksquare$

Then, by using the subadditivity of the matrix cumulant generating function \citep[see][Lemma 3.4]{tropp2012user}, we obtain the small-deviation bound for sums of random matrices:
\begin{theorem}\label{th:sumsLaplace}
Let $\{{\textbf X}_1,{\textbf X}_2,\cdots,{\textbf X}_K\}$ be a finite sequence of independent random Hermitian matrices. Then for any $\epsilon > 0$,
  \begin{equation}\label{eq:thsums}
   {\mathbb{P}}\left\{\lambda_{\max}\Big(\sum_k{{\textbf X}_k}\Big)\leq \epsilon\right\}\leq \inf_{\theta>0}\left\{{\rm e}^{\theta \epsilon}\cdot \exp\left(\lambda_{\max}\Big(\sum_k\log {\mathbb{E}}\,{\rm e}^{-\theta {\textbf X}_k}\Big)\right)\right\}.
  \end{equation}
  \end{theorem}
{\it Proof:}
By combining Lemma \ref{lem:laplace} and Lemma 3.4 of \citep{tropp2012user}, we have for any $\theta>0$,
 \begin{align*}
      {\mathbb{P}}\left\{\lambda_{\max}\Big(\sum_k{{\textbf X}_k}\Big)\leq \epsilon\right\}&\leq \frac{1}{d}\cdot {\rm e}^{\theta \epsilon}\cdot {\mathbb{E}}\,{\rm tr}\big(\,{\rm e}^{-\theta \sum_k{{\textbf X}_k}}\big)\\
       &\leq \frac{1}{d}\cdot {\rm e}^{\theta \epsilon}\cdot {\rm tr}\left(\exp\Big(\sum_k\log {\mathbb{E}}\,{\rm e}^{-\theta {\textbf X}_k}\Big)\right)\\
& \leq \frac{1}{d}\cdot {\rm e}^{\theta \epsilon} \cdot d \cdot \lambda_{\max}\left(\exp\Big(\sum_k\log {\mathbb{E}}\,{\rm e}^{-\theta {\textbf X}_k}\Big)\right)\\
& ={\rm e}^{\theta \epsilon}\cdot \exp\left(\lambda_{\max}\Big(\sum_k\log {\mathbb{E}}\,{\rm e}^{-\theta {\textbf X}_k}\Big)\right).
    \end{align*}
Taking the infimum over $\theta>0$ completes the proof.\hfill$\blacksquare$

Note that the above small-deviation bound is independent of the matrix dimension $d$, and thus it is applicable to the scenarios of high-dimensional and even infinite-dimensional matrices. In addition, we also derive the following small-deviation bounds for sums of random matrices.

\begin{corollary}\label{Co:g(theta)}
Let $\{{\textbf X}_1,{\textbf X}_2,\cdots,{\textbf X}_K\}$ be a sequence of independent random Hermitian matrices. Assume that there are a function $g(\theta)$ and a sequence $\{{\textbf A}_k\}$ of fixed Hermitian matrices such that
\begin{equation}\label{eq:mgf}
 {\mathbb{E}}\,{\rm e}^{-\theta {\textbf X}_k} \preceq {\rm e}^{g(\theta)\cdot {\textbf A}_k},\quad \forall\;\theta>0.
  \end{equation}
\begin{enumerate}[(i)]
\item Define the scalar parameter
  \begin{equation*}
    \eta_1:=\lambda_{\max}\Big(\sum_k{\textbf A}_k\Big).
  \end{equation*}
If $g(\theta)>0$, then for any $\epsilon > 0$,
  \begin{equation}\label{eq:cog(theta1)}
    {\mathbb{P}}\left\{\lambda_{\max}\Big(\sum_k{{\textbf X}_k}\Big)\leq \epsilon\right\}\leq \inf_{\theta>0}\big\{\exp \big(\theta\epsilon+g(\theta)\cdot\eta_1\big)\big\}.
  \end{equation}
\item Define the scalar parameter
  \begin{equation*}
    \eta_2:=\lambda_{\min}\Big(\sum_k{\textbf A}_k\Big).
  \end{equation*}
If $g(\theta)<0$, then for any $\epsilon > 0$,
  \begin{equation}\label{eq:cog(theta2)}
    {\mathbb{P}}\left\{\lambda_{\max}\Big(\sum_k{{\textbf X}_k}\Big)\leq \epsilon\right\}\leq \inf_{\theta>0}\big\{\exp \big(\theta\epsilon+g(\theta)\cdot\eta_2\big)\big\}.
  \end{equation}
\end{enumerate}
  \end{corollary}
{\it Proof:}
It follows from \eqref{eq:mgf} that
\begin{equation*}
  \log  {\mathbb{E}}\,{\rm e}^{-\theta {\textbf X}_k} \preceq {g(\theta)\cdot {\textbf A}_k},
\end{equation*}
and substituting it into Theorem \ref{th:sumsLaplace} leads to the result \eqref{eq:cog(theta1)}. Then, the fact $\lambda_{\max}(-{\textbf X})=-\lambda_{\min}({\textbf X})$ leads to the result \eqref{eq:cog(theta2)}. This completes the proof.\hfill$\blacksquare$

By using the logarithm operation, we then obtain another small-deviation bound for sums of random matrices:

\begin{corollary}\label{Co:sumsLaplace2}
Let $\{{\textbf X}_1,{\textbf X}_2,\cdots,{\textbf X}_K\}$ be a sequence of independent random Hermitian matrices. Then for any $\epsilon > 0$,
  \begin{equation*}\label{eq:cosums2}
    {\mathbb{P}}\left\{\lambda_{\max}\Big(\sum_k{{\textbf X}_k}\Big)\leq \epsilon\right\}\leq \inf_{\theta>0}\exp \left(\theta \epsilon+K\cdot\log\,\lambda_{\max}\Big(\frac{1}{n}\sum_{k=1}^K{\mathbb{E}}\,{\rm e}^{-\theta{\textbf X}_k}\Big)\right).
  \end{equation*}
  \end{corollary}
{\it Proof:}
Since the matrix logarithm is operator concave, for each $\theta>0$, we have    \begin{equation*}
      \sum_{i=1}^K\log\,{\mathbb{E}}\,{\rm e}^{-\theta{\textbf X}_k}=K \cdot \frac{1}{K}\sum_{i=1}^K\log\,{\mathbb{E}}\,{\rm e}^{-\theta{\textbf X}_k}\preceq K\cdot \log \Big(\frac{1}{K}\sum_{i=1}^K{\mathbb{E}}\,{\rm e}^{-\theta{\textbf X}_k}\Big).
    \end{equation*}
According to \eqref{eq:thsums}, we then arrive at
\begin{equation*}
    {\mathbb{P}}\left\{\lambda_{\max}\Big(\sum_k{{\textbf X}_k}\Big)\leq \epsilon\right\}\leq \frac{1}{d}\cdot {\rm e}^{\theta \epsilon}\cdot {\rm tr}\exp\left(K\cdot \log \Big(\frac{1}{K}\sum_{i=1}^K{\mathbb{E}}\,{\rm e}^{-\theta{\textbf X}_k}\Big)\right).
  \end{equation*}
Since the trace of a matrix can be bounded by $d$ times of its maximum eigenvalue, taking the infimum over $\theta>0$ completes the proof.\hfill$\blacksquare$

The following presents the relationship between one random psd matrix and a sum of psd random matrices.
\begin{lemma}\label{lem:sums}
Let $\{{\textbf X}_1,{\textbf X}_2,\cdots,{\textbf X}_K\}$ be a sequence of independent random Hermitian psd matrices. Then for any $\epsilon > 0$,
  \begin{equation*}\label{eq:one-sum}
    {\mathbb{P}}\left\{\lambda_{\max}\Big(\sum_k{\textbf X}_k\Big)\leq \epsilon\right\}\leq \prod_k{\mathbb{P}}\Big\{\lambda_{\max}({\textbf X}_k)\leq \epsilon\Big\}\leq {\mathbb{P}}\Big\{\lambda_{\max}({\textbf X}_k)\leq \epsilon\Big\}.
  \end{equation*}
  \end{lemma}
{\it Proof:}
 Since $\{{\textbf X}_1,{\textbf X}_2,\cdots,{\textbf X}_K\}$ are psd, we have
 \begin{align*}
 {\mathbb{P}}\left\{\lambda_{\max}\Big(\sum_k{{\textbf A}_k}\Big)\leq \epsilon\right\}&\leq {\mathbb{P}}\Big\{\max\Big(\lambda_{\max}({\textbf X}_1),\cdots,\lambda_{\max}({\textbf X}_K)\Big)\leq \epsilon\Big\}\nonumber\\
& =\prod_k{\mathbb{P}}\Big\{\lambda_{\max}({\textbf X}_k)\leq \epsilon\Big\}\nonumber\\
&\leq {\mathbb{P}}\Big\{\lambda_{\max}({\textbf X}_k)\leq \epsilon\Big\}.
 \end{align*}
The last inequality holds for any $k=1,2\cdots,K$. This completes the proof.\hfill$\blacksquare$

 This lemma shows that the small-deviation probability for sums of random matrices can be bounded by using the small-deviation probability for one single matrix. This fact suggests that the small-deviation bound could be independent of the size of matrix sequence, while this phenomenon will not arise in the large-deviation scenario.

\section{Small-deviation Inequalities for Positive Semi-Definite Random Matrices}\label{sec:psd}

In this section, we present several types of small-deviation inequalities for the largest eigenvalue of sums of independent random psd matrices. Similar to the scalar version of small-deviation inequalities, there remains a challenge to bound the term $\mathbb{ E}{\rm e}^{-\theta {\textbf X}_k}$. Here, we adapt some methods to handle this issue.

First, we introduce the negative moment estimate for the largest eigenvalue to derive a small-deviation inequality for sums of random matrices:
\begin{theorem}\label{th:sumsMoment}
Let $\{{\textbf X}_1,{\textbf X}_2,\cdots,{\textbf X}_K\}$ be a sequence of independent random Hermitian psd matrices. Given a $p>0$, if there exists a positive constant $C_p$ such that
 \begin{equation*}\label{eq:sumsmoment}
   \left[\lambda_{\max}\Big({\sum_k{{\mathbb{E}}\textbf X}_k}\Big)\right]^{-p}<C_p,
 \end{equation*}
then there holds that for any $\epsilon> 0$,
\begin{equation*}\label{eq:sums moment}
   {\mathbb{P}}\left\{\lambda_{\max}\Big({\sum_k{\textbf X}_k}\Big)\leq \epsilon\right\}\leq C_p\epsilon^p.
 \end{equation*}
 \end{theorem}
{\it Proof:}
It follows from Jensen's inequality that
   \begin{equation*}
    {\mathbb{E}}\,\left(\lambda_{\max}\Big({\sum_k{\textbf X}_k}\Big)\right)^{-p}\leq \left({\mathbb{E}}\,\lambda_{\max}\Big({\sum_k{\textbf X}_k}\Big)\right)^{-p}\leq  \left[\lambda_{\max}\Big({\sum_k{{\mathbb{E}}\textbf X}_k}\Big)\right]^{-p}<\infty.
   \end{equation*}
Then, the Markov's inequality yields
   \begin{align*}
      {\mathbb{P}}\left\{\lambda_{\max}\Big({\sum_k{\textbf X}_k}\Big)\leq \epsilon\right\}& ={\mathbb{P}}\left\{\left[\lambda_{\max}\Big({\sum_k{\textbf X}_k}\Big)\right]^{-p}\geq \epsilon^{-p}\right\}\nonumber\\
& \leq  \epsilon^p \cdot {\mathbb{E}}\left(\lambda_{\max}\big({\sum_k{\textbf X}_k}\big)\right)^{-p}\leq C_p\epsilon^p.
    \end{align*}
    This completes the proof.\hfill$\blacksquare$

In this theorem, we impose an assumption that the negative moment of $\lambda_{\max}\big({\sum_k{\textbf X}_k}\big)$ is bounded. In general, this assumption is mild and can be satisfied in most cases. The following small-deviation results are derived under that condition that the eigenvalues of the matrices $\{{\textbf X}_k\}$ are bounded:

\begin{theorem}\label{th:chernoff2}
Let $\{{\textbf X}_1,{\textbf X}_2,\cdots,{\textbf X}_K\}$ be a sequence of independent random Hermitian psd matrices such that $\lambda_{\max}({\textbf X}_k)\leq L$ ($\forall\; k=1,2,\cdots,K$) almost surely. Then for any $\epsilon > 0$,
 \begin{equation}\label{eq:chernoff2}
     {\mathbb{P}}\left\{\lambda_{\max}\Big({\sum_k{\textbf X}_k}\Big)\leq \epsilon\right\}\leq\left(\frac{\mu}{\epsilon}\right)^{{\epsilon}/{L}}\cdot \exp\left(\frac{\epsilon-\mu}{L}\right),
 \end{equation}
 where
 \begin{equation*}\label{eq:mu2}
  \mu:=\lambda_{\min}\left(\sum_k{\mathbb{E}}{\textbf X}_k\right).
  \end{equation*}
Furthermore, there holds that for any $\epsilon > 0$,
 \begin{equation}\label{eq:chernoff3}
     {\mathbb{P}}\left\{\lambda_{\max}\Big({\sum_k{\textbf X}_k}\Big)\leq \epsilon\right\}\leq \left(\frac{1}{\epsilon}\right)^{K\varepsilon/L}\cdot \left(\prod_{k=1}^K\mu_{k}\right)^{\epsilon/L}\cdot \exp\left(\frac{K\epsilon-\sum_k\mu_{k}}{L}\right),
  \end{equation}
where
\begin{equation*}\label{eq:mu k2}
  \mu_{k}=\lambda_{\min}({\mathbb{E}}{\textbf X}_k).
 \end{equation*}
   \end{theorem}
{\it Proof:}
For any $\theta>0$ and $x\in [0,L]$, there holds that
   \begin{equation*}
     {\rm e}^{-\theta x}\leq 1+\frac{{\rm e}^{-\theta L}-1}{L}\cdot x\leq \exp\Big(\frac{{\rm e}^{-\theta L}-1}{L}\cdot x\Big).
   \end{equation*}
According to transfer rule, we have,
     \begin{equation}\label{eq:cgf2}
       \log{\mathbb{E}}\,{\rm e}^{-\theta {\textbf X}_k}\preceq \frac{{\rm e}^{-\theta L}-1}{L}{\mathbb{E}}{\textbf X}_k.
     \end{equation}
By substituting \eqref{eq:cgf2} into the Corollary \ref{Co:g(theta)}, we then have for any $\theta>0$,
     \begin{align*}
      {\mathbb{P}}\left\{\lambda_{\max}\Big(\sum_k{{\textbf X}_k}\Big)\leq \epsilon\right\}&\leq {\rm e}^{\theta \epsilon}\cdot \exp\left(\lambda_{\max}\Big(\sum_k\frac{{\rm e}^{-\theta L}-1}{L}{\mathbb{E}}{\textbf X}_k\Big)\right)\nonumber\\
& =\exp\left(\theta \epsilon+\lambda_{\max}\Big(\frac{{\rm e}^{-\theta L}-1}{L}\sum_k{\mathbb{E}}{\textbf X}_k\Big)\right)\nonumber\\
& =\exp\left(\theta \epsilon+\frac{{\rm e}^{-\theta L}-1}{L}\cdot \lambda_{\min}\Big(\sum_k{\mathbb{E}}{\textbf X}_k\Big)\right)\nonumber\\
& =\exp\left(\theta \epsilon+\frac{{\rm e}^{-\theta L}-1}{L}\cdot \mu\right).
     \end{align*}
The infimum is achieved at $\theta=\frac{1}{L}\log({\frac{\mu}{\epsilon}})$, which leads to the result of \eqref{eq:chernoff2}.

Moreover, the combination of Lemma \ref{lem:laplace} and \eqref{eq:cgf2} leads to
  \begin{equation*}
      {\mathbb{P}}\left\{\lambda_{\max}({{\textbf X}_k})\leq \epsilon\right\} \leq \left(\frac{\mu_{k}}{\epsilon}\right)^{{\epsilon}/{L}}\cdot \exp\left(\frac{\epsilon-\mu_{k}}{L}\right).
     \end{equation*}
Then, the result \eqref{eq:chernoff3} is derived from Lemma \ref{lem:sums}. This completes the proof.\hfill$\blacksquare$

Actually, the above results are derived from the geometric point of view, where the term ${\rm e}^{-\theta x}$ is bounded by the linear function $1+\frac{{\rm e}^{-\theta L}-1}{L}$ for any $x\in[0,L]$. Finally, we study the small-deviation inequalities for random matrix series $\sum_k x_k{\textbf A}_k$, which is a sum of fixed Hermitian psd matrices ${\textbf A}_k$ weighted by random variables $x_k$.

\begin{theorem}\label{th:poly}
Let $\{{\textbf A}_1,{\textbf A}_2,\cdots,{\textbf A}_K\}$ be a sequence of fixed Hermitian psd matrices, and $\{x_1,x_2,\cdots,x_K\}$ be a finite sequence of independent variables. If there exist the constants $C>0$ and $\alpha>0$ such that
\begin{equation}\label{eq:poly}
  {\mathbb{E}}\,{\rm e}^{-\theta x_k}\leq C\cdot \theta^{-\alpha},
\end{equation}
then there holds that for any $0<  \epsilon< \frac{K\alpha}{{\rm e}}\cdot\sqrt[\alpha]{\frac{K }{C\nu}}$,
\begin{equation}\label{eq:bound poly}
  {\mathbb{P}}\left\{\lambda_{\max}\Big({\sum_k x_k{\textbf A}_k}\Big)\leq \epsilon\right\}\leq \left(\frac{{\rm e} \epsilon}{K \alpha}  \right)^{\alpha K} \cdot\left(\frac{C\nu}{K}\right)^K,
\end{equation}
where
\begin{equation*}\label{eq:nu}
 \nu=\lambda_{\max}\left(\sum_k{{\textbf A}_k}^{-\alpha}\right).
\end{equation*}
Furthermore, for any $0<\epsilon <( \alpha / {\rm e}) \cdot C^{-\frac{1}{\alpha}}\cdot \left(\prod_{k=1}^K\nu_k\right)^{-\frac{1}{\alpha K}} $,
\begin{equation}\label{eq:bound poly2}
  {\mathbb{P}}\left\{\lambda_{\max}\Big({\sum_k x_k{\textbf A}_k}\Big)\leq \epsilon\right\}\leq \left(\prod_{k=1}^K\nu_k\right)\cdot C^K\cdot \left(\frac{{\rm e}\epsilon}{\alpha}\right)^{K \alpha},
\end{equation}
where
\begin{equation*}\label{eq:nu_k}
  \nu_k=\lambda_{\max}({{\textbf A}_k}^{-\alpha}).
\end{equation*}
\end{theorem}
{\it Proof:}
According to transfer rule, we have,
  \begin{equation}\label{eq:poly mgf}
    {\mathbb{E}}\,{\rm e}^{-\theta x_k{\textbf A}_k}\preceq C\cdot (\theta {\textbf A}_k)^{-\alpha}.
  \end{equation}
By substituting \eqref{eq:poly mgf} into the Corollary \ref{Co:sumsLaplace2}, we then have for any $\theta>0$,
  \begin{align*}
      {\mathbb{P}}\left\{\lambda_{\max}\Big(\sum_kx_k{{\textbf A}_k}\Big)\leq \epsilon\right\}&\leq \exp \left(\theta \epsilon+K\cdot\log\,\lambda_{\max}\Big(\frac{1}{K}\sum_{k=1}^KC\cdot (\theta {\textbf A}_k)^{-\alpha}\Big)\right)\nonumber\\
& =\exp \left(\theta \epsilon+K\cdot\log\,\lambda_{\max}\Big(\frac{C\theta^{-\alpha}}{K}\sum_{k=1}^K {\textbf A}_k^{-\alpha}\Big)\right)\nonumber\\
& =\exp \left(\theta \epsilon+K\cdot\log\frac{C\nu}{K\theta^{\alpha}}\right).
     \end{align*}
The infimum will be attained at $\theta=\frac{\alpha K}{\epsilon}$, and it leads to the result of \eqref{eq:bound poly}. Moreover, the combination of Lemma \ref{lem:laplace} and \eqref{eq:poly mgf} leads to
  \begin{equation*}
      {\mathbb{P}}\left\{\lambda_{\max}(x_k{{\textbf A}_k})\leq \epsilon\right\} \leq C\cdot \nu_k \cdot \left(\frac{{\rm e}\epsilon}{\alpha}\right)^\alpha.
     \end{equation*}
Then, the result \eqref{eq:bound poly2} is resulted from Lemma \ref{lem:sums}. This completes the proof.\hfill$\blacksquare$

The above results hold under the condition \eqref{eq:poly} that ${\mathbb{E}}\,{\rm e}^{-\theta x_k}$ has a power-type upper bound of $ C\cdot \theta^{-\alpha}$. This condition is mild and we refer to \citet{liten} for the details. Moreover, to keep the results \eqref{eq:bound poly} and \eqref{eq:bound poly2} non-trivial, their right-hand sides should be less than {\it one}, and thus we arrive at $\epsilon< \frac{K\alpha}{{\rm e}}\cdot\sqrt[\alpha]{\frac{K }{C\nu}}$ and $\epsilon <( \alpha / {\rm e}) \cdot C^{-\frac{1}{\alpha}}\cdot \left(\prod_{k=1}^K\nu_k\right)^{-\frac{1}{\alpha K}}$, respectively.

\section{Conclusion}

In this paper, we present the small-deviation inequalities for the largest eigenvalues of sums of random matrices. In particular, we first give some basic results on small-deviation inequalities for random matrices.
We then study the small-deviation inequalities for sums of independent random psd matrices. Different from the large-deviation inequalities for random matrices, our results are independent of the matrix dimension $d$ and thus can be applicable to the scenarios of high-dimensional and even infinite-dimensional matrices. In addition, by using the Hermitian dilation \cite[see][Section 2.6]{tropp2012user}, our small-deviation results can also be extended to the scenario of non-Hermitian random matrices.



\section*{Acknowledgement}

This work is partially supported by the Fundamental Research Funds for the Central Universities: DUT13RC(3)068 and DUT17LK46; the National Natural Science Foundation of China: 11401076 and 61473328; Dalian High Level Talent Innovation Support Program: 2015R057.

%


\end{document}